\documentclass[a4paper,12pt,reqno]{amsart}
\usepackage{amssymb}
\usepackage{graphicx}
\usepackage{amsmath}
\usepackage{color}
\usepackage{enumitem}
   \nonstopmode \numberwithin{equation}{section}
\newtheorem{thm}{Theorem}[section]
 \newtheorem{cor}[thm]{Corollary}
 \newtheorem{lem}[thm]{Lemma}
 \newtheorem{prop}[thm]{Proposition}
 \theoremstyle{definition}
 \newtheorem{defn}[thm]{Definition}
 \theoremstyle{remark}
 \newtheorem{rem}[thm]{Remark}
 \newtheorem{rems}[thm]{Remarks}
  \newtheorem{ques}[thm]{Question}
 \newtheorem{ex}[thm]{Example}

 \newtheorem{prob}[thm]{Problem}
 \numberwithin{equation}{section}

\def\be{\begin{equation}}
\def\ee{\end{equation}}

\newcommand{\blem}{\begin{lem}}
\newcommand{\elem}{\end{lem}}
\newcommand{\bthm}{\begin{thm}}
\newcommand{\ethm}{\end{thm}}
\newcommand{\bcor}{\begin{cor}}
\newcommand{\ecor}{\end{cor}}
\newcommand{\beg}{\begin{ex}}
\newcommand{\eeg}{\end{ex}}
\newcommand{\bprop}{\begin{prop}}
\newcommand{\eprop}{\end{prop}}
\newcommand{\bdefn}{\begin{defn}}
\newcommand{\edefn}{\end{defn}}
\newcommand{\bprob}{\begin{prob}}
\newcommand{\eprob}{\end{prob}}
\newcommand{\bpf}{\begin{proof}}
\newcommand{\epf}{\end{proof}}

\newcommand{\brem}{\begin{rem}}
\newcommand{\erem}{\end{rem}}
\newcommand{\brems}{\begin{rems}}
\newcommand{\erems}{\end{rems}}

\newcommand{\beq}{\begin{eqnarray}}
\newcommand{\beqq}{\begin{eqnarray*}}
\newcommand{\eeq}{\end{eqnarray}}
\newcommand{\eeqq}{\end{eqnarray*}}
\newcommand{\bqn}{\begin{ques}}
\newcommand{\eqn}{\end{ques}}

\newcommand{\vp}{\varphi}

\newcommand{\C}{{\mathbb C}}
\newcommand{\A}{{\mathcal A}}
\newcommand{\D}{{\mathbb D}}

\newcommand{\N}{{\mathbb N}}

\newcounter{minutes}
\divide\time by 60
\newcounter{hours}
\multiply\time by 60 \addtocounter{minutes}{-\time}
\begin{document}
\title
{ Multiplicative operators on analytic function spaces}

\author{Kanha Behera}
\address{Kanha Behera, Department of Mathematics and Statistics, Indian Institute of Technology, Kanpur - 208016, India.}
\email{beherakanha2560@gmail.com}

\author{Junming Liu}
\address{Junming Liu, School of Mathematics and Statistics, Guangdong University of Technology, Guangzhou, Guangdong - 510520, People's Republic of China.}
\email{jmliu@gdut.edu.cn}

\author{P. Muthukumar}
\address{P. Muthukumar,
Department of Mathematics and Statistics, Indian Institute of Technology, Kanpur - 208016, India.}
\email{pmuthumaths@gmail.com, muthu@iitk.ac.in}

%


\begin{abstract}
 H. J. Schwartz proved in his thesis (1969) that a nonzero bounded operator on Hardy spaces $(H^p, 1\leq p\leq\infty)$ is almost  multiplicative if and only if it is a composition operator. But, his proof has a gap.  In this article, we show that his result is not correct for $H^\infty$ and we fill the gap for $H^p, 1\leq p<\infty.$ Further, we prove that on several classical spaces such as  the Bloch space, the little Bloch space, Besov spaces $B_p$ for $p>1$, and weighted Bergman spaces an  operator is almost multiplicative if and only if it is a composition operator.  Finally, we give a  complete characterization of those composition operators that are multiplicative with respect to the Duhamel product of analytic functions.
\end{abstract}

\thanks{
File:~\jobname .tex,
          printed: \number\day-\number\month-\number\year,
          \thehours.\ifnum\theminutes<10{0}\fi\theminutes
}
\maketitle
\tableofcontents

\pagestyle{myheadings}
\markboth{KANHA BEHERA, JUNMING LIU and P. MUTHUKUMAR}{Multiplicative operators on analytic function spaces}

\section{Introduction}\label{Introprelim}
Let $\D=\{z\in\mathbb{C} : |z|<1\}$ denote the open unit disc in the complex plane and let
 $\mathbb{T}=\{z\in\mathbb{C} : |z|=1\}$ be its boundary. Let $\mathcal{A}$ be a space of analytic functions on $\D$.
 For an analytic self-map $\varphi$ of $\D$, the associated composition operator $C_\varphi$ on $\mathcal{A}$ is defined by
  \[
    C_\varphi(f)=f\circ\varphi,\qquad f\in \mathcal{A}.
  \]
%
  A bounded linear operator $T$ on $\mathcal{A}$ is called \emph{almost multiplicative} if
                   \[T(fg)=T(f) T(g)\]
 whenever $f,g,fg\in\mathcal{A}$. If $\A$ is an algebra under pointwise multiplication, then an almost
  multiplicative operator is simply a multiplicative operator. It is easy to see that composition
   operators on $\A$ are always almost multiplicative.

 In the late 1930s, the relation between multiplicative bounded linear functionals on a Banach algebra
  and its maximal ideals was studied by I. Gelfand, and this relation formed the basis of the field
   of commutative algebra (see \cite[Chapter 7]{ConwayFnAnaBook}). Later investigations into
   multiplicative linear operators on normed algebras were also carried out (see \cite{multopalgIKap}
   and related publications). In another direction, closely related $\delta$-multiplicative linear
   functional and $\delta$-multiplicative linear operators on different normed algebras were also
    studied (see \cite{AE}, \cite{JK1}, \cite{JK2}, \cite{JB}, \cite{HuaLin}, \cite{PSemrel} and references therein).

 Study of almost multiplicative operators on normed spaces, which are not algebras was initiated by
 H. J. Schwartz in his thesis \cite{Schwartzthesis} in 1969. He proved that on all Hardy spaces a
 nonzero bounded linear operator is almost multiplicative if and only if it is a composition operator.
 In Section $\ref{Schwartzthesis}$, we show that Schwartz's result for  $H^\infty$
 is not  correct and his proof for the Hardy spaces $H^p$ spaces, $1\leq p<\infty$, has a
 small gap which we fill along with an alternative proof.

 In Section $\ref{almultbesovsec}$, we give a characterization of almost multiplicative operators on the
 Besov spaces $B_p$ for $1<p<\infty$ and see that they are all composition operators. In
 Section $\ref{blochsec}$, we get the same result for almost multiplicative operators on the little
 Bloch space and the Bloch space and in Section $\ref{bergmansec}$, we get the same result for weighted
  Bergman spaces.

 In Section $\ref{multopcompsec}$, we consider the class of analytic function spaces of form $\psi\A$
 where $\A$ is an algebraically consistent function space and $\psi$ is a suitable analytic function
 such that $\psi\A\subset\A$. We prove that any nonzero bounded linear almost multiplicative operator
  on such $\psi\A$ is a composition operator.

 In Section $\ref{duhalgsec}$, we consider the Duhamel product of two analytic functions under which
  many well-known analytic function spaces form an algebra called Duhamel algebra. Here we completely
   characterize the composition operators on Duhamel algebras that are almost multiplicative operators
   with respect to the Duhamel product.

To end this section, let us give a general remark on almost multiplicative operators on normed spaces.
 \brem\label{remT1eq1}
Let $\A$ be a normed space of analytic functions on $\D$ with constant function $\textbf{1}\in \A$ and
$T:\A\to \A$ be a nonzero almost multiplicative operator. Then $T\textbf{1}=\textbf{1}$ as
$$T\textbf{1} = T\textbf{1}\cdot T\textbf{1} \implies T\textbf{1}(\textbf{1}-T\textbf{1})= 0$$
and this gives us $T\textbf{1}=\textbf{1}$, because if $T\textbf{1}=0$ then $T\equiv 0$.
\erem

\section{Revisiting Schwartz's result}\label{Schwartzthesis}

Recall that for $0< p<\infty$, $H^p$ is the space containing analytic functions on $\D$ satisfying
$$
\|f\|_{p}:=\sup\limits_{ r \in [0,1)}\left \{\frac{1}{2\pi}\int_0^{2\pi}
|f(re^{i\theta})|^{p}\,d\theta\right \}^{\frac{1}{p}}<\infty
$$
and $H^\infty$ is the space of all bounded analytic functions on $\D$, i.e.,
$$
\|f\|_{\infty}:=\sup\limits_{z\in\D}|f(z)|<\infty \text{   for all } f\in H^\infty.
 $$
For $1\leq p \leq \infty$, $H^p$ is a Banach space.
 By a consequence of Littlewood's subordination theorem, it is well-known
  that composition operators map $H^p$ spaces into themselves and thus composition operators are
   bounded on $H^p$ (see \cite[Corollary, Page-29]{DurenHpspaces} or \cite[Corollary 3.7]{Cowenbook}).

Schwartz proved in his thesis, in \cite[Theorem 1.3]{Schwartzthesis}, that composition operators
 are the only multiplicative operators on $H^\infty$.  But, this is incorrect.

\bthm\label{Hinfnotcompmult}
 Schwartz's result \cite[Theorem 1.3]{Schwartzthesis} is not correct. That is, there are multiplicative
 operators on $H^\infty$ which are not composition operators.
\ethm
\bpf
Fix  $w\in\mathbb{T}$. Then, the function $z-w\in H^\infty$ is a non-invertible element as
$\lim_{z\to w}(z-w)=0$. Therefore, the scalar $w$ is in the spectrum of the function $z$ in $H^\infty$ and by
an elementary result in commutative algebra (see \cite[Chapter VII, Theorem 8.6]{ConwayFnAnaBook})
 there exists a multiplicative bounded linear operator (functional)
 $T:H^\infty\to \C \cong \C\cdot\textbf{1}\subseteq H^\infty$ such that $T(z-w\cdot\textbf{1})=0$
 i.e., $Tz=w$. It gives that $T$ is a multiplicative bounded linear operator on
  $H^\infty$, which is not a composition
 operator. Otherwise, 
  \[
  w=Tz=C_\vp (z)=\vp
 \] and thus $w\in\D$, which is not possible.
\epf

\brem
Indeed, in \cite[Theorem 1.3]{Schwartzthesis}, Schwartz actually proved the following:
 Let $T:H^\infty\to H^\infty$ be a nonzero bounded linear multiplicative operator. Then $\vp:=Tz$ maps
$\D$ to $\overline{\D}$. Moreover, if $\vp$ is a self map of $\D$ then $T=C_\vp$ on $H^\infty$.
\erem
An analogous result for the disc algebra $A(\D)$, the closure
  of the polynomials in $H^\infty$, is given below. For sake of completeness,
we give its proof.

\bthm\label{multHinfty}
Let $T:A(\D)\to A(\D)$ be a nonzero bounded linear multiplicative operator. Then $\vp:=Tz$ maps
$\D$ to $\overline{\D}$. Further,
\begin{itemize}
  \item[(i)] If $\vp:\D\to\D$, then $T=C_\vp$ on $A(\D)$.
  \item[(ii)] If $\vp$ is a uni-modular constant $c$, then $Tf =f(c)$ on  $A(\D)$.
\end{itemize}
\ethm
\bpf
Suppose $T$ is a multiplicative operator on $A(\D)$. By Remark \ref{remT1eq1}, we have $T\textbf{1}=\textbf{1}$. Let $\vp=Tz$.
Since $T$ is almost multiplicative , we have
$$T(z^n)=(Tz)^n=\vp^n \ \ \text{for all} \ n\in\N.$$
As $\|z^n\|_\infty=1$, therefore
$$\|\vp^n\|_\infty\leq \|T\|  \ \ \text{for all} \ n\in\N.$$
If $|\vp(z)|=\delta$ for some $z\in\D$, then $\delta^n \leq \|\vp^n\|_\infty\leq \|T\|$, for all $n$.
 We can not have $\delta >1$, because in that case $\delta^n\to \infty$ leading to a contradiction.
 Hence, $\vp$ maps $\D$ into $\overline{\D}$ and therefore because of maximum modulus principle we have
 two cases, either $\vp$ is a self-map of $\D$ or it is a uni-modular constant.

Let us consider the first case, that is, $\vp:\D\to\D$. Fix $w\in\D$.
For $f\in A(\D)$, we have
 $$f-f\circ\vp(w)=(z-\vp(w))g,$$
 where $g\in A(\D)$. By
 applying $T$ on both sides of above equation, we get
 $$ Tf-f\circ\vp(w)=(\vp-\vp(w))Tg.$$
  Right hand side vanishes at $w$ and therefore $Tf(w)=f\circ\vp(w)$. Since $f\in A(\D)$ and $w\in\D$ are arbitrarily chosen,
  we have $T=C_\vp$ on $A(\D)$.

  Now, if $\vp\equiv c$, an uni-modular constant, then by linearity and multiplicativity of $T$
  we have $Tp=p(c)$ for any polynomial $p$. Since the disc algebra, $A(\D)$, is the closure of
  all polynomials, we have $Tf=f(c)$ on $A(\D)$.
\epf
In his thesis \cite[page 14]{Schwartzthesis}, H. J. Schwartz proved that for $1\leq p<\infty$ a
bounded operator $T$ on $H^p$ is almost multiplicative if and only if it is a composition operator.
It comes to our notice that Schwartz's proof has a small but not so obvious to ignore gap in it.

Suppose $T$ is a nonzero almost multiplicative bounded operator on $H^p$ and let $Tz=\vp$. Schwartz has proved that $|\vp(z)|\leq 1$
 for all $z \in\D$. Maximum modulus principle yields that $\vp$ maps $\D$ to itself or it is a
 uni-modular constant function. Schwartz has not addressed the case:  $\vp$ is a uni-modular constant.

We use the following result to fill the gap in the proof of Schwartz.
\bprop\label{Log1_zc}
For any scalar $\lambda\in \mathbb{T}$, the map $f(z)=\log\frac{1}{1-\lambda z}\in H^p$ for all
 $p\in (0,\infty)$. Further, for $1\leq p<\infty$, Taylor's polynomials
  $f_n(z)=\sum_{k=1}^{n}\frac{(\lambda z)^k}{k}$ of $f$ converges to $f$ in $H^p$ norm.
\eprop
\bpf
Observe that $(\log(\frac{1}{1-z}))'=\frac{1}{1-z}\in H^p$ for all $p\in (0,1)$. Therefore,
by \cite[Theorem 5.12]{DurenHpspaces}, we have $\log(\frac{1}{1-z})\in H^q$ for all
$q=p/(1-p)\in(0,\infty)$. Since composition operators maps $H^q$ to itself for all
$q\in (0,\infty)$, we have
 $$
 f(z)=\log\left(\frac{1}{1-\lambda z}\right)=
 \sum_{k=1}^{\infty}\frac{(\lambda z)^k}{k}\in H^p  \text{ for all } p\in (0,\infty).
 $$

Now, by \cite[Corollary 3]{ZhuNormCvg}, $f_n(z) = \sum_{k=1}^{n}\frac{(\lambda z)^k}{k}$
converges to $\log(\frac{1}{1-\lambda z})$ in $H^p$ norm for $1<p<\infty$.
 As the inclusion map $H^2$ to $H^1$ is bounded, $f_n$ converges to $f$ in $H^1$ norm.
\epf
Now, we not only fill the gap in Schwartz's proof  but also give an alternative proof of it.
\bthm
Fix $1\leq p<\infty$. Let $T$ be a nonzero bounded linear operator on $H^p$. Then, $T$ is
 almost multiplicative if and only if $T$ is a composition operator.
\ethm
\bpf
Composition operators are always almost multiplicative. For the other way, suppose $T$ is a nonzero  almost multiplicative operator on $H^p$.
 Define $\vp:=Tz$ as earlier, so that
$$
Tz^n=\vp^n \ \ \text{for all} \ n\in\N.
$$
As $\vp=Tz\in H^p$, it is well-known that the radial limit of $\vp$
satisfies $\vp(e^{i\theta})\in L^p[0,2\pi]$. \\
\textbf{Claim:} $|\vp(e^{i\theta})|\leq 1$ for almost every $\theta\in [0,2\pi]$.

For a contradiction, suppose $|\vp(e^{i\theta})|>1$ on $E\subset[0,2\pi]$, such that the
 Lebesgue measure of $E$, $m(E)>0$. Then for some $\delta>1$, it is possible to find
  $E'\subset E$ with $m(E')>0$ such that $|\vp(e^{i\theta})|>\delta$ on $E'$.
  For all $n\in\N$, one has
$$
\|T\|^p \geq \frac{\|Tz^n\|_p}{\|z^n\|_p}=
\|\vp^n\|_p^p \geq \frac{1}{2\pi}\int_{E'}|\vp(e^{i\theta})|^{np}
d\theta \geq \delta^{np}\frac{m(E')}{2\pi}.
$$
This is impossible since $\delta^{np}\to\infty$ as $n$ increases. Hence we must
 have $|\vp(e^{i\theta})|\leq 1$ almost everywhere on $[0,2\pi]$. As
 $\vp\in H^p$ and $\vp(e^{i\theta})\in L^\infty(\partial\D)$, we get
   $\vp\in H^\infty$  (by \cite[Theorem 2.11]{DurenHpspaces}) with
  $\|\vp\|_\infty=\|\vp\|_{L^\infty(\partial\D)}\leq 1$, that is, $\vp:\D\to\overline{\D}$.
  Now there are two cases, either $\vp:\D\to\D$ or $\vp$ is a uni-modular constant function.

Suppose $\vp\equiv c$, where $|c|=1$. Since $T$ is linear and almost multiplicative, we have
$Th=h(c)$ for any polynomial $h$. Let $f(z)=\log\frac{1}{1-\overline{c} z}$ and
$f_n(z)=\sum_{k=1}^{n}\frac{(\overline{c} z)^k}{k}$ be Taylor's polynomials of $f$.
By Proposition \ref{Log1_zc}, we have $f_n$ converges to $f$ in $H^p$ norm. Now, since
$T$ is a bounded linear operator on $H^p$, therefore $Tf_n= f_n(c)$ converge to $Tf$ in $H^p$ norm.
But $\|Tf_n\|_p=|f_n(c)|=\sum_{k=1}^{n}\frac{1}{k}$ converges to $\infty$, which is a
 contradiction. Hence $\vp$ cannot be a uni-modular constant. And hence $\vp$ is
 a self-map of $\D$.

 For the later part of  the proof, Schwartz has used the denseness of polynomials in
  $H^p$ spaces, we give an alternative proof. Fix $w\in\D$. For a given $f\in H^p$,
  there exist $g\in H^p$ such that
 $$f-f\circ\vp(w)=(z-\vp(w))g.$$
 Existence of $g$ is follows by Riesz factorization theorem
 \cite[Theorem 2.5]{DurenHpspaces}. 
  Almost multiplicativity of  $T$ gives
 $$ Tf-f\circ\vp(w)=(\vp-\vp(w))Tg.$$
  Consequently, $Tf(w)=f\circ\vp(w)$. Since $f\in H^p$ and $w\in\D$ are arbitrarily
  chosen, we get $T=C_\vp$ on $H^p$.
\epf

\section{Multiplicative operators on $B_p$}\label{almultbesovsec}
For $p\in (1,\infty)$, the Besov space $B_p$ consists of all analytic functions on $\D$ that satisfies
 $$\int_{\D}(1-|z|^2)^{p-2}|f'(z)|^pdA < \infty$$
 with the norm
 $$\|f\|_{B_p} = |f(0)|+\left(\int_{\D}(1-|z|^2)^{p-2}|f'(z)|^pdA\right)^{1/p},$$
 where $dA=\frac{1}{\pi}dxdy$ is the normalized area measure on $\D$. For $p=1$, the
 space $B_1$ is the set of all analytic functions on $\D$ that satisfies
 $$\int_{\D}|f''(z)|dA<\infty,$$
 the norm here is
 $$\|f\|_{B_1}=|f(0)|+|f'(0)|+\int_{\D}|f''(z)|dA.$$
Before proving the main result of this section, let us give a general lemma for arbitrary
analytic function spaces. The lemma considers some sufficient conditions under which
composition operators are the only almost multiplicative operators on analytic function spaces.
\blem\label{almultgenfnsp}
Let $X$ be a Banach space of analytic functions on $\D$ and let $T:X\to X$
 be a nonzero bounded linear operator. Assume that:
\begin{enumerate}
\item Polynomials are dense in $X$,
\item $Tz$ maps $\D$ into itself,
\item Evaluation maps $K_a:X\to\C$ defined by $K_{a}f:= f(a)$, where
$a\in\D$, are bounded.
\end{enumerate}
Then, $T$ is almost multiplicative if and only if $T$ is a composition operator on $X$.
\elem
\bpf
Since composition operators are always almost multiplicative, it is enough to prove the other way. Let $\vp=Tz$
as before. Almost multiplicativity
 of $T$ gives that $Tp=p\circ\vp$ for any polynomial $p$. Let $f\in X$. By
 denseness of polynomials in $X$, we have a sequence of polynomials $p_n$ converging
 to $f$. Boundedness of $T$ implies that $p_n\circ\vp$ converge to $Tf$. 
 Since $\vp$ is a self-map and evaluation maps are bounded, applying $K_{\vp(z)}$
  and $K_z$ we get that $p_n\circ\vp(z)$ converges to $f\circ\vp(z)$ as
   well as $Tf(z)$ for arbitrary $z\in\D$. Hence, we $T=C_\vp$ on $X$.
\epf
\brem\label{blochfnineqrem}
For $f\in\mathcal{B}$, by \cite[Theorem 5.1.6]{ZhuOpThFnSpBook}, we have for any $a\in\D$
$$|f(a)-f(0)|\leq \frac{1}{2}\log\frac{1+|a|}{1-|a|}\|f\|_\mathcal{B}.$$
Consequently, evaluation maps are bounded on $\mathcal{B}$. As $B_p\subset \mathcal{B}$
for all $1<p<\infty$, the inclusion maps are bounded by the closed graph theorem. 
Thus for all $f\in B_p$, we have $\|f\|_\mathcal{B}\leq C_p\|f\|_{B_p}$, where $C_p$ is a positive constant.
Therefore, evaluation maps are also bounded on $B_p$. 
As a result, one has that norm convergence 
 in Besov space implies pointwise convergence.
\erem
\bthm\label{almultbesov}
Let $1< p <\infty$ and $T: B_p\to B_p$ be a nonzero bounded linear operator. Then, $T$ is
 almost multiplicative  if and only if $T$ is a composition operator.
\ethm
\bpf
Suppose $T$ is almost multiplicative. 
Let $\vp=Tz$.
 Thus, we get
$$T\left(\frac{z^n}{n}\right)=\frac{\vp^n}{n} \ \ \text{for all } n\in\N.$$
Now,
\[
\begin{split}
\|z^n/n\|_{B_p}^p & =
\frac{1}{\pi}\int_{\D}(1-|z|^2)^{p-2}\,|z|^{(n-1)p}dA\\
& =  \int_{0}^{1}r^{(n-1)p}(1-r^2)^{p-2}2rdr 
\\
& =B\left(\frac{(n-1)p}{2}+1,p-1\right),
\end{split}
\]
where $B(\cdot,\cdot)$ is the beta function. Recall that, if $x$ is large and $y$ is fixed we have
the approximate formula for the beta function, $B(x,y)\sim \Gamma (y)x^{-y}.$ This relation
 and $p-1>0$ implies that $\|z^n/n\|_{B_p}\to 0$ as $n\to\infty$, in particular,
  $\|z^n/n\|_{B_p}\leq M$ (a constant) for all $n\in\N$. Now,
$$\|\vp^n/n\|_{B_p}=\|Tz^n/n\|_{B_p}\leq \|z^n/n\|_{B_p} \|T\| \ \leq M\|T\|,
\ \ \text{for all} \ n\in\N.$$
\textbf{Claim:} $|\vp(z)|< 1$ on $\D$.

\textbf{Case 1:} Suppose $\vp\equiv c$, a constant. Now,
$$\dfrac{|c|^n}{n} =\|\vp^n/n\|_{B_p} 
\leq M\|T\|, \ \ \text{for all} \ n\in\N.$$
If $|c|>1$, then $\frac{|c|^n}{n}\to\infty$ as $n$ increases, which leads to a contradiction.
Next assume $|c|=1$. At first, let us suppose $c=1$. Take $f(z)=\left( \log\frac{2}{1-z}\right)^\gamma$,
 for some $0<\gamma < 1-1/p$, then $f\in B_p$ (see \cite[Theorem 1]{besovfneg}).
 Let $p_n$ be the $n$th Taylor polynomial of $f$. By \cite[Corollary 6]{ZhuNormCvg}, $p_n$
 converges to $f$ in Besov space norm. Thus, $Tp_n=p_n(1)$ converges. If $p_n(1)$
 converges then by Abel's limit theorem (see \cite[Theorem 2.5]{AhlforsComplexAna}),
 radial limit of $f$ at $1$ exists, which is a contradiction. Since $B_p$ is
  M\"{o}bius-invariant, for any $|c|=1$ we will get a similar contradiction.
Therefore, $|c|<1$.

\textbf{Case 2:} Assume that $\vp$ is  non constant.
Suppose that $|\vp(a)|>1$ for some $a\in\D$. By continuity of $\vp$, choose a
neighbourhood $B\subseteq \D$ of $a$ such that $|\vp|>\delta >1$ on $B$. Thus,
\[
\begin{split}
\left\|\vp^n/n\right\|_{B_p}^p & 
= \int_{\D}|\vp|^{(n-1)p}|\vp'(z)|^p(1-|z|^2)^{p-2} dA \\
& \geq \int_{B}|\vp|^{(n-1)p}|\vp'(z)|^p(1-|z|^2)^{p-2} dA
\geq \delta^{(n-1)p}K,
\end{split}
\]
where $K=\int_{B}|\vp'(z)|^p(1-|z|^2)^{p-2} dA$. Note that $K>0$ because,
$B$ is of nonzero Lebesgue measure and $\vp$ is nonconstant leads
$\vp'$ is zero only at countably many points in $\D$. Now,
$$\delta^{n-1}K^{1/p}\leq \|\vp^n/n\|_{B_p}\leq M\|T\|, \ \ \text{for all} \ n\in\N,$$
which is not possible as the left side goes to $\infty$ as $n\to\infty$. Hence, we must
 have $|\vp(z)|\leq 1$. Since $Tz=\vp$ is non-constant, it is a self-map of $\D$.
Since polynomials are dense in $B_p$ and
 evaluation maps are bounded, Lemma \ref{almultgenfnsp} completes the proof.
\epf

\section{Multiplicative operators on $\mathcal{B}$}\label{blochsec}

The Bloch space, $\mathcal{B}$, is a Banach space of analytic functions on $\D$ such that
$$\|f\|_{\mathcal{B}}= |f(0)| + \sup_{z\in\D}(1-|z|^2)|f'(z)| < \infty .$$
 Also, $\|f\|_{\mathcal{B}} \leq |f(0)|+ \|f\|_{\infty}$ and therefore
 $H^\infty \subset \mathcal{B}$ (see \cite[Proposition 5.1.2]{ZhuOpThFnSpBook}).
The little Bloch space, $\mathcal{B}_0$, is the closed subspace of $\mathcal{B}$
which consists of functions such that $(1-|z|^2)f'(z)\to 0$ as $|z|\to 1^{-}$.
$\mathcal{B}_0$ is the closure of polynomials in $\mathcal{B}$. For any self-map
 $\vp$, the composition operator $C_\vp$ is bounded on $\mathcal{B}$. Also,  $C_\vp$
 is bounded on $\mathcal{B}_0$ if and only if $\vp\in\mathcal{B}_0$
 (see \cite[Page 2]{compopblochsp}).

\bthm\label{almultlilbloch}
Let $T: \mathcal{B}_0\to \mathcal{B}_0$ be a nonzero bounded linear operator. Then,
$T$ is almost multiplicative  if and only if $T=C_\vp$ for some self-map
$\vp\in\mathcal{B}_0$.
\ethm
\bpf
If $\vp\in\mathcal{B}_0$, then $C_\vp$ is bounded on $\mathcal{B}_0$ and
composition operators always almost multiplicative. Conversely, suppose $T$ is an
 almost multiplicative bounded linear operator on $\mathcal{B}_0$.
 Take $\vp=Tz$ as before. Then,
$$Tz^n=\vp^n \ \ \text{for all} \ n\in\N.$$
By the relation between norms of $\mathcal{B}$ and $H^\infty$ (given above),
one has
$$\|z^n\|_{\mathcal{B}} \leq \|z^n\|_{\infty}= 1\ \ \ \text{for all } n\in\N .$$

\textbf{Claim:} $|\vp(z)|< 1$ on $\D$.

\textbf{Case 1:} Suppose, $\vp$ is a non constant analytic function.
Now,
$$\|\vp^n\|_{\mathcal{B}}= \|Tz^n\|_{\mathcal{B}}\leq
 \|z^n\|_{\mathcal{B}}\, \|T\| \ \leq \|T\| \ \ \text{for all} \ n\in\N,$$
and thus for all $n>1$ and for any $z\in\D$, we have
\beq\label{eqbdvp}
(1-|z|^2)|\vp(z)|^{n-1}\,|\vp'(z)|\leq \frac{1}{n}\|T\|.
\eeq
If $\vp'(z)\neq 0$ and $|\vp(z)|>1$ for some $z\in\D$, then the left side of equation
 \eqref{eqbdvp} tends to infinity as $n$ increases, which is a contradiction. Also, since
 $\vp$ is non constant and analytic on $\D$, $\vp'$ can have at most countably many zeros
  in $\D$ otherwise identity theorem forces $\vp'$ to be zero leading to a contradiction.
   Therefore $|\vp(z)|\leq 1$ for all but countably many points in $\D$ and continuity of
    $\vp$ yields $|\vp|\leq 1$ on $\D$. Now, by maximus modulus principle $\vp$ is either
     a self-map or a uni-modular constant. Since the later can not be true $\vp$ is a
     non-constant self-map in this case.

\textbf{Case 2:} Suppose, $\vp\equiv c$, a constant. Now,
$$|c|^n =\|\vp^n\|_{\mathcal{B}} \leq \|z^n\|_{\mathcal{B}} \|T\|\leq \|T\|, \ \
 \text{for all} \ n\in\N.$$
If $|c|>1$, then the left side goes to infinity as $n$ increases, which is a contradiction.

Suppose, $|c|=1$. For $n\in\N$, consider the following polynomial
$$p_n(z)=\sum_{k=0}^{n-1}\frac{(\overline{c}z)^{2k+1}}{2k+1}.$$
Then,
$$p_n'(z)=\sum_{k=0}^{n-1}(\overline{c}z)^{2k}$$
and
$$|p_n'(z)| \leq \sum_{k=0}^{n-1}|z|^{2k}=\frac{1-|z|^{2n}}{1-|z|^2}.$$
Therefore,
$$\|p_n\|_{\mathcal{B}} = |p_n(0)| +\sup_{z\in\D}(1-|z|^2)|p_n'(z)| \leq 1$$
for all $n\in\N$. But we have,
$T(p_n) = p_n(c) = \sum_{k=0}^{n-1}\frac{1}{2k+1}$ diverging to $\infty$ as $n$ increases.
This is a contradiction to the fact that $T$ is a bounded linear operator. Hence $|c|<1$.

Therefore, in any case, $\vp$ maps $\D$ into itself and $z\in\mathcal{B}_0$
gives $Tz=\vp\in\mathcal{B}_0$. As a consequence $C_\vp:\mathcal{B}_0\to \mathcal{B}_0$
 is a bounded linear operator. And as $Tz^n = C_\vp z^n$ for all $n\in\N\cup\{0\}$
 therefore $T=C_\vp$ on polynomials. Finally, since polynomials are dense in
 $\mathcal{B}_0$, we have $T=C_\vp$ on $\mathcal{B}_0$.
\epf
Before we proceed to the almost multiplicative operators on $\mathcal{B}$,
we need following results.
\blem\label{blochfnineq}
If $f\in\mathcal{B}$, then $\sup\limits_{z\in\D}(1-|z|^2)|f(z)|<\infty$.
\elem
\bpf
By Remark \ref{blochfnineqrem}, we have for $f\in\mathcal{B}$ and $z\in\D$,
$$(1-|z|^2)|f(z)|\leq \frac{(1-|z|^2)}{2}\log\frac{1+|z|}{1-|z|}\|f\|_\mathcal{B}+(1-|z|^2)|f(0)|.$$
To prove our claim we only need to show
$$\sup_{z\in\D}(1-|z|)\log\frac{1}{1-|z|}<\infty.$$
And the above is true as
\begin{align*}
\lim_{r\to 1^-}(1-r)\log\frac{1}{1-r}=\lim_{x\to\infty}\frac{\log x}{x}=\lim_{x\to\infty}\frac{1}{x}=0.
\end{align*}
The desired result follows.
\epf

 \blem\label{blochfactor}
 Suppose $f\in\mathcal{B}$ and $w\in\D$, then $f(z)-f(w)=(z-w)g(z)$ for some $g\in\mathcal{B}$.
 \elem
 \bpf
 For $f\in\mathcal{B}$ and $w\in\D$, let us define
 $$
  g(z):=\left\{
\begin{array}{cl}
\dfrac{f(z)-f(w)}{z-w}& \mbox{ if } z\neq w,\\
f'(w) & \mbox{ if } z=w.
\end{array}
\right.
  $$

 Since $f(z)-f(w)$ is analytic in $\D$ and has a zero at $w$, therefore $g$ is also
 analytic in $\D$. We claim that $\sup\limits_{z\in\D}(1-|z|^2)|g'(z)|<\infty$.
 Since $g'$ is bounded in a closed ball with center $w$ contained in $\D$, let us
 consider $|z-w|>\delta$ for some $\delta >0$. Now,
 \begin{align*}
 (1-|z|^2)|g'(z)|&\leq \frac{1-|z|^2}{\delta^2}|(z-w)f'(z)-(f(z)-f(w))| \\
                       &\leq \frac{2}{\delta^2}(2(1-|z|^2)|f'(z)|+(1-|z|^2)|f(z)|+2|f(w)|),
 \end{align*}
 and as $f\in\mathcal{B}$, Lemma \ref{blochfnineq} ensures that $\sup_{z\in\D}
 (1-|z^2|)|g'(z)|<\infty$.
It completes the proof.
 \epf

 \bthm\label{almultbloch}
Let $T: \mathcal{B}\to \mathcal{B}$ be a nonzero bounded linear operator. Then,
$T$ is almost multiplicative if and only if  $T$ is a composition operator.
\ethm
\bpf
Proceeding the same way as in Theorem \ref{almultlilbloch} we get $Tz=\vp\in\mathcal{B}$
is a self-map and $T=C_\vp$ on polynomials. For $f\in \mathcal{B}$ and $w\in\D$,
by Lemma \ref{blochfactor} we have
 $$f-f\circ\vp(w)=(z-\vp(w))g,$$
 where $g\in\mathcal{B}$. Now, since $T$ is almost multiplicative, we have
 $$ Tf-f\circ\vp(w)=(\vp-\vp(w))Tg.$$
  Therefore, $Tf(w)=f\circ\vp(w)$ by evaluating at $w$ both sides. Since $f\in \mathcal{B}$ and $w\in\D$ are
  arbitrarily chosen, we get $T=C_\vp$ on $\mathcal{B}$.
\epf

\section{Multiplicative operators on $A^p_\alpha$}\label{bergmansec}
For $p>0$ and $\alpha >-1$, we define the weighted Bergman space $A^p_\alpha$ to be the space of
 analytic functions on $\D$ that satisfies,
$$ \|f\|^p_{A^p_\alpha}:= \int_{\D}|f(z)|^p(\alpha +1)(1-|z|^2)^\alpha dA < \infty ,$$
where $dA$ is the normalized area measure on $\D$. For $1\leq p <\infty$ and $\alpha >-1$,
it is known that $A^p_\alpha$'s are Banach spaces and polynomials are dense in them.

In \cite[Theorem 3.2]{Aalpha2AGupta} it is proved that any almost multiplicative operator
 on $A^2_\alpha$, for $\alpha>-1$, is a composition operator. In the following theorem, we
 provide an alternative and simpler proof of their result, and at the same time
 generalize it to $A^p_\alpha$ for $1\leq p <\infty$ and $\alpha >-1$.

\bthm\label{bergalmult}
Let $1\leq p <\infty$ and $\alpha>-1$. 
Then, a nonzero bounded linear operator on $A^p_\alpha$ is almost multiplicative
if and only if it is a composition operator.
\ethm

\bpf
Composition operators $f\to f\circ\vp$ are always almost multiplicative. For the
converse part, suppose $T$ is almost multiplicative. 
 Take $\vp=Tz$. Thus,
$$Tz^n=(Tz)^n=\vp^n \ \ \text{for all} \ n\in\N.$$
\[
\begin{split}
\|z^n\|_{A^p_\alpha}^p & = 
\frac{1}{\pi}\int_{\D}|z^n|^p (\alpha +1)(1-|z|^2)^\alpha dA \\
& = 2(\alpha + 1)\int_{0}^{1}r^{np+1}(1-r^2)^\alpha dr
\\
& = (\alpha +1)B\left(\frac{np}{2}+1,\alpha +1\right),
\end{split}
\]
where $B(\cdot,\cdot)$ is the beta function. Now if $x$ is large and $y$ is fixed we have the  approximate formula for the beta function, $B(x,y)\sim \Gamma (y)x^{-y}.$ This relation
 and $\alpha + 1>0$ implies that $\|z^n\|_{A^p_\alpha}\to 0$ as $n\to\infty$, in particular,
 $(\|z^n\|_{A^p_\alpha})$ is a bounded sequence, say bounded by a constant
 $M>0$.
Thus,
$$\|\vp^n\|_{A^p_\alpha}\leq \|z^n\|_{A^p_\alpha}\, \|T\|\leq M\|T\|, \ \ \text{for all} \ n\in\N.$$

\textbf{Claim:} $|\vp(z)|< 1$ on $\D$.

Suppose that $|\vp(a)|>1$ for some $a\in\D$. By continuity of $\vp$, choose a neighbourhood $B\subseteq \D$ of
$a$ such that $|\vp|>\delta >1$ on $B$. Thus,
$$\|\vp^n\|_{A^p_\alpha}^p=\int_{\D}|\vp|^{np}(\alpha +1)(1-|z|^2)^\alpha dA
\geq \int_{B}|\vp|^{np}(\alpha +1)(1-|z|^2)^\alpha dA\geq \delta^{np}K,$$
where $K=\int_{B}(\alpha +1)(1-|z|^2)^\alpha dA> 0$.
 Now,
$$\delta^{n}K^{1/p}\leq \|\vp^n\|_{A^p_\alpha}\leq M\|T\|, \ \
 \text{for all} \ n\in\N,$$
which is not possible as the left side goes to $\infty$ as $n\to\infty$.
 Hence, we must have $|\vp(z)|\leq 1$.

But if $|\vp(a)|=1$ for some $a\in\D$, then by maximum modulus principle $\vp$
is a uni-modular constant. In this case,
$$K_1=\|\vp^n\|_{A^p_\alpha}= \|Tz^n\|_{A^p_\alpha}\leq \|z^n\|_{A^p_\alpha} \|T\|
, \ \ \text{for all} \ n\in\N,$$
where $\|\vp^n\|_{A^p_\alpha}=K_1=\int_{\D}(\alpha +1)(1-|z|^2)^\alpha dA> 0$ for all $n$.
This is not possible because $\|z^n\|_{A^p_\alpha}$ and hence the right side goes to $0$ as
 $n\to\infty$, whereas the left side is a fixed positive constant for all $n$.

Therefore, $\vp$ maps $\D$ into itself. As a consequence $C_\vp:A^p_\alpha\to A^p_\alpha$ is
a bounded linear operator. And as $Tz^n = C_\vp z^n$ for all $n\in\N\cup\{0\}$ therefore
 $T=C_\vp$ on polynomials. Finally, denseness of polynomials in $A^p_\alpha$ gives
 the required result,
  $T=C_\vp$ on $A^p_\alpha$.
\epf


\section{General setting}\label{multopcompsec}

The following theorem is a generalization of \cite[Theorem 1.4]{Cowenbook} and it gives
 us an alternate method for characterization of almost multiplicative operators on analytic function spaces.

 \bthm\label{Aconjugcompop}
 Suppose $\mathcal{A}$ be a normed space of analytic functions on $\D$ such that
 \begin{itemize}
   \item[(i)] all evaluation
 maps, $K_z:f\to f(z)$ for $z\in\D$, are bounded on $\mathcal{A}$,
   \item[(ii)] there exists $g\in\mathcal{A}$ and a univalent function
    $\psi$ analytic on $\D$ such that $\psi g\in\mathcal{A}$. 
 \end{itemize} Then, a bounded linear operator $T$ on $\mathcal{A}$
  is a composition operator if and only if the set $\{K_z:z\in\D\setminus S\}$
    is mapped into evaluation maps by $T^*$
 for some countable $S\subset\D$, such that none of the limit points of $S$ are in $\D$.
 \ethm
 \bpf
 If $T$ is a composition operator, that is $T=C_\vp$ for some self-map $\vp$ on $\D$,
 then for any $x\in\D$ and $f\in\A$ we have,
 $$(T^*K_x)(f)=K_x(Tf)=f\circ\vp(x)=K_{\vp(x)}(f),$$
 that is, $T^*K_x=K_{\vp(x)}$ for any $x\in\D$. Hence, the set of all evaluation maps,
 $\{K_z:z\in\D\}$, is mapped into evaluation maps by $T^*$.

 Now, for the converse part, assume that $T^*$ maps the collection $\{K_z:z\in\D\setminus S\}$ into $\{K_z:z\in\D\}$,
  where $S\subset\D$ is countable and $\overline{S}\cap\D=S$.
  Thus, $\D\setminus S$ is an open set.
   Define $\vp:\D\setminus S \to \D$ by $T^*K_z=K_{\vp(z)}$.

 \textbf{Claim:} $\vp$ has analytic extension upto $\D$.

 For any $x\in\D\setminus S$ and $f\in\A$, we have
 $$(Tf)(x)=K_x(Tf)=(T^*K_x)(f)=K_{\vp(x)}(f)=f\circ\vp(x).$$
 As $Tf\in\A$ is analytic on $\D$, we have $f\circ\vp$ is analytic on the open set
  $\D\setminus S$ for any $f\in\A$.

 Also, given that there exists $g\in\A$ and an univalent function $\psi$ such that $\psi g\in\A$.
  By the previous logic, $g\circ\vp$ and $(\psi\circ\vp)(g\circ\vp)$ are analytic on
   $\D\setminus S$. Therefore, $\psi\circ\vp$ is analytic on $\D\setminus S$ as all
   the zeros of $g\circ\vp$ are also the zeros of $(\psi\circ\vp)(g\circ\vp)$ of same
   or greater multiplicities. Now, as $\psi$ is an univalent analytic function on $\D\setminus S$
    therefore $\psi^{-1}$ is also an univalent analytic function from $\psi(\D\setminus S)$ onto
     $\D\setminus S$. Hence, $\vp=\psi^{-1}\circ(\psi\circ\vp)$ is analytic on $\D\setminus S$.
     As $\vp$ maps $\D\setminus S$ into $\D$, points of $S$ can have only removable
     singularities of $\vp$. Therefore, $\vp$ can be extended to $\D$, so 
      that $\vp$ is analytic on $\D$.

 Now as $f\circ\vp$ is also analytic on $\D$, by identity theorem we have
 $$Tf=f\circ\vp ~\text{for all } f\in\A$$
 that is, $T=C_\vp$ on $\A$.
 \epf

 \bdefn\label{algconsdefn}
  Let $\A$ be a Banach space of analytic functions on $\D$ with all evaluation
   maps are bounded. We say that $\A$ is \textbf{algebraically consistent} if for every
   nonzero almost multiplicative bounded linear functional on $\A$ is given by evaluation
   at some point of $\D$.
 \edefn

 \brem\label{fnspalgcon}
 Since it is proved in Sections \ref{Schwartzthesis}, \ref{almultbesovsec}, \ref{blochsec}
  and \ref{bergmansec} that Hardy spaces, Besov spaces, the little Bloch space, the Bloch
   space and weighted Bergman spaces have the property that a nonzero bounded linear
   operator is almost multiplicative if and only if it is a composition operator and
   also  these spaces have constant functions, it is easy to see that these spaces
   are algebraically consistent.
 \erem

 A multiplication operator on a analytic function space $\A$ is defined by
$$M_h(f)(z)=h(z)f(z),$$
where $h$  is an analytic function  on $\D$. 
If all evaluation  maps are bounded on $\A$, then 
it is easy to see that multiplication operators
are closed operators and hence by the closed graph theorem, $M_h$ on space $\A$ is bounded
 if and only if $M_h(\A)\subset\A$ i.e., $h$ is a multiplier of $\A$. Also it is easy to
 see that polynomials are multipliers of Hardy spaces and weighted Bergman spaces.
  We can apply Lemma \ref{blochfnineq} to see that polynomials are multipliers of
  the Bloch space and the little Bloch space. Now, let us give a lemma that will help
   us in showing that polynomials are multipliers of Besov spaces.

\blem\label{Besovfnineq}
If $f\in B_p$, for $1<p<\infty$, then $\int_{\D}(1-|z|^2)^{p-2}|f(z)|^pdA < \infty$.
\elem
\bpf
If $f\in B_p$, for $1<p<\infty$, then by \cite[Theorem 9]{besovsppapzhu} we have
the growth estimate,
$$|f(z)|\leq C\|f\|_{B_p}\left(\log\frac{2}{1-|z|^2}\right)^{1-1/p},$$
where $C$ is a constant. Therefore, we have,
$$\int_{\D}(1-|z|^2)^{p-2}|f(z)|^pdA \leq C\|f\|_{B_p}\int_{\D}(1-|z|^2)^{p-2}
\left(\log\frac{2}{1-|z|^2}\right)^{p-1}dA.$$
Changing the integral on the right to polar form, it becomes
$$
2\int_{0}^{1}(1-r^2)^{p-2}\left(\log\frac{2}{1-r^2}\right)^{p-1}rdr.
$$
Now, by substituting $\log\frac{2}{1-r^2} = x$ we will see that the above improper
 integral is equivalent to
$$2^{p-1}\int_{\log2}^{\infty}x^{p-1}e^{-(p-1)x}dx$$
and after a routine substitution it is easy to see that this integral is bounded by
 a constant times $\Gamma(p)$.
\epf
As an immediate consequence, we get the following result.
\bthm\label{polymultBesov}
Suppose $g\in H^\infty$ such that $g'\in H^\infty$, then for any $f\in B_p$,
 where $1<p<\infty$, we have $fg\in B_p$.
\ethm
\bpf
Note that
\[
\begin{split}
\int_{\D}(1-|z|^2)^{p-2}|(fg)'(z)|^pdA  & 
\leq 2^p \int_{\D}(1-|z|^2)^{p-2}|(fg')(z)|^pdA +\\
&~~~ 2^p \int_{\D}(1-|z|^2)^{p-2}|(f'g)(z)|^pdA.
\end{split}
\]
The boundedness of the integrals in the RHS of the above inequality follows from
Lemma \ref{Besovfnineq} and the facts that $g,g'\in H^\infty$ and $f\in B_p.$
\epf
\brem\label{rempolymultBes}
As an application of Theorem \ref{polymultBesov}, we have polynomials are multipliers
of $B_p$, where $1<p<\infty$.
\erem

Now we give a new class of function spaces which are algebraically consistent and in
turn they have the property that any nonzero bounded linear operator on them is almost
 multiplicative if and only if it is a composition operator. In the  following theorem,
 we need polynomials to be dense in the Banach space $\A$. Polynomials are dense in Hardy
 spaces $H^p$ for $1< p<\infty$, weighted Bergman
 spaces, Besov spaces $B_p$ for $1<p<\infty$ (see \cite[Corollary 3, 4, 6]{ZhuNormCvg}) and the little Bloch space (by definition).
  Also as discussed before, these spaces are algebraically consistent Banach spaces and
   polynomials are multipliers of these spaces. Hence, we can give the following result
   where $\A$ can be any one of the above mentioned spaces as a particular case.

 \bprop\label{subspalgcons}
 Suppose $\A$ be an algebraically consistent functional Banach space of analytic functions
  on $\D$ such that polynomials are dense in $\A$ and suppose $M_p(\A)\subset\A$ for all
  polynomials $p$. If $\psi\A\subset\A$ then  $\psi\A$ is also an algebraically
  consistent normed space.
 \eprop
 \bpf
 Since evaluation maps are bounded on $\A$ and $\psi\A\subset\A$, evaluation
 maps are also bounded on $\psi\A$. Suppose $k$ be a nonzero, bounded linear functional
 on $\psi\A$ such that $k(fg)=k(f)k(g)$ whenever $f$, $g$, and $fg$ are in $\psi\A$.

Since the set of polynomials, P, is dense in $\A$ and $M_\psi$ is bounded on $\A$,
we have $\psi P$ is dense in $\psi\A$. There exists $p\in P$ such that $k(\psi p)\neq 0$,
 otherwise $k\equiv 0$ on $\psi P$ and therefore on $\psi\A$ too. As $p\A\subset\A$,
 we have   $\psi p\A\subset p\A\subset\A$. Thus standard application  of closed graph theorem
  yields that the multiplication operator
$M_{\psi p}$ is bounded operator on $\A$. That is, there exist $C>0$ with
\[
\|\psi pf\|_{\A}\leq C \|f\|_{\A} \mbox{~ for all~} f\in \A.
\]

 Define $\tilde{k}$ on $\A$ by
 $$\tilde{k}(f):=\frac{k(\psi pf)}{k(\psi p)}.$$
  Now, $\tilde{k}$ is a bounded linear functional as
  $$|\tilde{k}(f)|=\left|\frac{k(\psi pf)}{k(\psi p)}\right|\leq
  \frac{\|k\|}{|k(\psi p)|}\|\psi pf\|_{\A}\leq C_1 \|f\|_{\A},$$
  and $\tilde{k}|_{\psi\A}=k$ as
  $$\tilde{k}(\psi f)=\frac{k(\psi p\psi f)}{k(\psi p)}=
  \frac{k(\psi p)k(\psi f)}{k(\psi p)}=k(\psi f).$$
 Also whenever $f$, $g$, and $fg$ are in $\A$ we have
  $$\tilde{k}(fg)=\frac{k(\psi pfg)k(\psi p)}{k(\psi p)^2}=
  \frac{k(\psi pf\psi pg)}{k(\psi p)^2}=\frac{k(\psi pf)}
  {k(\psi p)}\frac{k(\psi pg)}{k(\psi p)}=\tilde{k}(f)\tilde{k}(g).$$
  Therefore, $\tilde k$ is a nonzero, bounded linear functional on $\A$
   that is also almost multiplicative and since $\A$ is algebraically consistent, $\tilde{k}$
    and thus $k$ is a evaluation map on  $\psi\A$.
 \epf

 \bcor\label{almultopsubsp}
 Take $\A$ and $\psi$ as in Proposition \ref{subspalgcons}. Suppose $T$ be a nonzero,
  bounded linear operator on $\psi \A$. Then, $T$ is almost
  multiplicative if and only if $T$ is a composition operator.
 \ecor
 \bpf
 Suppose $T$ be a nonzero, almost multiplicative, bounded linear operator on $\psi\A$.
 Choose $f\in\A$
  such that  $T(\psi f)\not\equiv 0$. Hence by the identity theorem, $Z(T(\psi f))$, the set
  of all zeros of $T(\psi f)$ in $\D$, is a countable set without a limit point in $\D$.
   For $x\in\D\setminus Z(T(\psi f))$,
 $$T^*K_x(\psi f)=K_x(T(\psi f))=T(\psi f)(x)\neq 0,$$
 that is, for each $x\in\D\setminus Z(T(\psi f))$ the bounded linear functional $T^*K_x\not\equiv 0$.
  Also whenever $f,g,fg$ are in $\psi A$,
 $$T^*K_x(fg)=K_x(T(fg))=K_x(TfTg)=K_x(Tf)K_x(Tg)=T^*K_x(f)T^*K_x(g),$$
 that is $T^*K_x$ is almost multiplicative. As $\psi\A$ is algebraically consistent by
 Proposition \ref{subspalgcons}, $T^*K_x$ is a evaluation map on $\psi\A$ for any
 $x\in\D\setminus Z(T(\psi f))$. Hence, by Theorem \ref{Aconjugcompop}, $T$ is a
 composition operator on $\psi \mathcal{A}$. The converse part trivially follows.
 \epf
 \brem\label{remonCowen}
 Take $\A$ and $\psi$ as in Proposition \ref{subspalgcons}.
 Note that if $\psi$ has more than one zero in $\D$ then $\psi\A$ is not a functional Banach space
  by \cite[Definition 1.1]{Cowenbook}, as if $x,y$ are two distinct zeros of $\psi$ in 
  $\D$ then
   $f(x)=f(y)=0$ for all $f\in\psi\A$ but $x\neq y$. Even, in this also, we still have $T$ is almost
   multiplicative if and only if $T$ is a composition operator on $\psi\A$.
 \erem

 \section{Duhamel multiplicative operators}\label{duhalgsec}

 The Duhamel product $f\circledast g$ of analytic functions $f$ and $g$ on $\D$ is defined as,
 $$(f\circledast g)(z):= \frac{d}{dz}\int_{0}^{z}f(z-t)g(t)dt=\int_{0}^{z}f'(z-t)g(t)dt+f(0)g(z).$$

It is well-known that Hardy spaces \cite[Theorem 1]{DuhamelHp}, Bergman spaces
\cite[Theorem 2.3]{DuhamelBergman}, Wiener disc algebra \cite[Theorem 4]{DuhamelWiener},
$Q_p$ spaces and Morrey spaces \cite[Theorem 1]{DuhamelQpMorreyspaces}, Besov spaces
 \cite[Theorem 1.1]{DuhamelBesov} are commutative Banach algebras under the Duhamel
 product $\circledast$. Let $X$ be an analytic function space that is a Banach algebra
 under the Duhamel product, if an operator $T$ is multiplicative on $X$ with respect
 to the Duhamel product then we will say that $T$ is Duhamel multiplicative on $X$,
 i.e., $T(f\circledast g) = Tf\circledast Tg$  for all $f,g\in X$.

 \bprop\label{phi00}
 Suppose $\vp$ be a self-map of $\D$ and $X$ be a Banach algebra under the Duhamel product
 and the identity map $z\in X$. If $C_\vp$ is Duhamel multiplicative on $X$ then $\vp(0)=0$.
 \eprop
 \bpf
 We have $C_\vp(f\circledast g) = C_\vp f\circledast C_\vp g$  for all $f,g\in X$. In particular,
  by taking $f=g=z$, we get $\frac{\vp^2}{2}=\vp\circledast\vp$ that is,
 $$\frac{\vp^2(z)}{2}=\int_{0}^{z}\vp'(z-t)\vp(t)dt+\vp(0)\vp(z), \ \text{for all}~ z\in\D.$$
 By taking $z=0$,
 $\frac{\vp^2(0)}{2}=\vp^2(0)$, which is possible only when $\vp(0)=0$.
 \epf
 \brem
 The converse of Proposition \ref{phi00} is not true. For example,
  take $\vp=z^2$, so that $\vp(0)=0$ and $f=g=z$. But,
  $$ C_\vp(f\circledast g) = \frac{z^4}{2} \neq \frac{z^4}{6}=C_\vp f\circledast C_\vp g. $$
 \erem
 Now we give a class of self-maps $\vp$ for which $C_\vp$ is Duhamel multiplicative.
 \bprop\label{phiaz}
Let $X$ be a Banach algebra under the Duhamel product. If $\vp(z)=az,|a|\leq 1$ and
if $C_\vp$ is a bounded operator on $X$, then $C_\vp$ is Duhamel multiplicative on $X$.
 \eprop
 \bpf
 Let $f,g\in X$. Then,
 $$C_\vp f\circledast C_\vp g(z)= (f\circ\vp)\circledast (g\circ\vp)(z)=
  \int_{0}^{z}af'(a(z-t))g(at)dt+f(0)g(az).$$
 By a change of variable $s=at$, we get
 $$C_\vp f\circledast C_\vp g(z) =\int_{0}^{az}f'(az-s)g(s)ds +f(0)g(az)=
 f\circledast g(az).$$
 Therefore, $ C_\vp f\circledast C_\vp g = C_\vp(f\circledast g)$ for all $f,g\in X$.
 \epf
 The following result completely characterizes which composition operators are  Duhamel
  multiplicative.
 \bthm\label{phipoly}
 Let $X$ be a Banach algebra under the Duhamel product with the identity function $z\in X$ and
 $C_\vp$ is a bounded operator on $X$. Then, $C_\vp$ is Duhamel multiplicative on $X$
 if and only if $\vp(z)=az$ for some $|a|\leq 1$.
 \ethm
 \bpf
 If $\vp(z)=az$ for some $|a|\leq 1$, then $C_\vp$ is Duhamel multiplicative by
 Proposition \ref{phiaz}.

 For the converse part assume that $C_\vp$ is Duhamel multiplicative on $X$. Then,
  first we see that $\vp(0)=0$
 by Proposition \ref{phi00}. Suppose $\vp(z)=\sum_{n=1}^{\infty}a_nz^n$ in $\D$.
 We shall use induction to show that $a_j=0$ for $j\geq 2$.

  Let $f=g=z$. We have for each $z\in\D$, $(f\circledast g)\circ\vp(z)=(f\circ\vp)
  \circledast(g\circ\vp)(z)$, which implies that $$\vp(z)^2/2=\int_{0}^{z}\vp'(z-t)\vp(t)dt.$$
   By taking $\vp(z)=\sum_{n=1}^{\infty}a_nz^n$, we obtain
 \begin{equation}\label{eqpolyDuhmult}
   \frac{1}{2}\left(\sum_{n=1}^{\infty}a_nz^n\right)^2  =  \int_{0}^{z}\left(\sum_{n=1}^{\infty}na_n(z-t)^{n-1}\right)
   \left(\sum_{n=1}^{\infty}a_nt^n\right)dt.
 \end{equation}
%
By comparing $z^3$ terms on both sides, we get
\[a_1a_2 =
  \frac{2}{3}a_1a_2,\] which is possible only if $a_1a_2=0$. If $a_1\neq 0$, then $a_2$ must be 0.
  Now suppose that $a_1=0$. 
Now, by comparing $z^4$ terms on both sides, we see that  
\[  
\frac{1}{2}a_2^2 = \frac{1}{6}a_2^2,
\] which is possible only if $a_2=0$.
  Hence in any case we have $a_2=0$.

Assume that  $a_2=a_3=\cdots=a_k=0$ for $k\geq 2$. Now equation (\ref{eqpolyDuhmult}) becomes
\begin{equation}\label{eqpolyDuhmult1}
   \frac{1}{2}\left(a_1z + \sum_{n=k+1}^{\infty}a_nz^n\right)^2  =
   \int_{0}^{z}\left(a_1 + \sum_{n=k+1}^{\infty}na_n(z-t)^{n-1}\right)
   \left(a_1t + \sum_{n=k+1}^{\infty}a_nt^n\right)dt
 \end{equation}
 In the equation \eqref{eqpolyDuhmult1}, the $z^{k+2}$ term on the right side is 
 \[
\begin{split}
 = & \int_{0}^{z}\left( a_1a_{k+1}t^{k+1} + (k+1)a_{k+1}(z-t)^k\cdot a_1t\right) dt
\\
= & a_1a_{k+1}\frac{z^{k+2}}{k+2} + a_1a_{k+1}(k+1)\frac{1! \ k!}{(k+2)!}z^{k+2}
\\
=& \frac{2}{k+2}a_1a_{k+1}z^{k+2}.
\end{split}
\]
 
Here we have used the following formula, using change of parameter $t=sz$, for 
positive integers $m, n$, we get
\[
\begin{split}
 \int_{0}^{z}(z-t)^{m}\cdot t^ndt & = z^{m+n}\int_{0}^{z}(1-t/z)^{m}(t/z)^ndt
\\
& = z^{m+n+1}\int_{0}^{1}(1-s)^{m}s^nds
\\
& = z^{m+n+1}\, B(n+1,m+1)
\\
& = z^{m+n+1}\frac{n! \ m!}{(m+n+1)!}.
\end{split}
\]
Now, equating the coefficient of $z^{k+2}$ on the left and the right of equation \eqref{eqpolyDuhmult1}
 we get, 
 \[a_1a_{k+1}=\frac{2}{k+2}a_1a_{k+1},\]
  which is possible only if $a_1a_{k+1}=0$,
  because $\frac{2}{k+2}=1$ if and only if  $k=0$. Thus, either $a_1=0$ or $a_{k+1}=0$.
  Suppose $a_1=0$, then the $z^{2(k+1)}$ term in the right is
\begin{align*}
\int_{0}^{z}(k+1)a_{k+1}(z-t)^k\cdot a_{k+1}t^{k+1}dt = \frac{(k+1)k!(k+1)!}{(2k+2)!}a_{k+1}^2z^{2(k+1)}.
\end{align*}
Now, comparing coefficients of $z^{2(k+1)}$ on both side of the equation (\ref{eqpolyDuhmult1}),
we obtain
$$\frac{1}{2}a_{k+1}^2 = \frac{(k+1)!(k+1)!}{(2k+2)!}a_{k+1}^2.$$
Using induction principle, it is easy to see that $\frac{(k+1)!(k+1)!}{(2k+2)!}<\frac{1}{2}$ for all $k>1$.
Hence, we must have $a_{k+1}=0$. As an application of mathematical induction, 
we get 
\[
a_n=0, n\geq 2.
\] 
Therefore, $\vp(z)=a_1z$. As $\vp$ is a self-map of $\D$, we also have $|a_1|\leq1$.
 \epf
As all polynomials belongs to Hardy spaces, Bergman spaces, Wiener algebra, $Q_p$ spaces, Morrey
spaces and Besov spaces and also these spaces are Banach algebras under the Duhamel product the above
 result holds in these spaces.

The Bloch space, $\mathcal{B}$ is a Banach space of analytic functions on $\D$ such that
$$\|f\|_{\mathcal{B}}= |f(0)| + \sup_{z\in\D}(1-|z|^2)|f'(z)| < \infty .$$
  For $p\in (1,\infty)$, $Q_p=\mathcal{B}$, the Bloch space (see \cite[Cor 1.2.1]{Qclassesbook}).
  So the Bloch space is also a Duhamel algebra, but we give an alternate and direct proof of it.

 \bthm
 The Bloch space, $\mathcal{B}$  is a Banach algebra under the Duhamel product.
 \ethm
 \bpf
 Let $f,g\in\mathcal{B}$. Then,
 $$(f\circledast g)(z) =\int_{0}^{z}f'(z-t)g(t)dt+f(0)g(z) = \int_{0}^{z}g(z-t)f'(t)dt+f(0)g(z).$$
Consequently,
\beq\label{eqduhderiv}
(f\circledast g)'(z) = \int_{0}^{z}g'(z-t)f'(t)dt+g(0)f'(z)+f(0)g'(z).
\eeq
Recall that for any $h\in\mathcal{B}$, we have $|h(0)|\leq\|h\|_{\mathcal{B}}$ and  $(1-|z|^2)|h'(z)|\leq\|h\|_{\mathcal{B}}$ for all
 $z\in\D$. Considering the polar form,
we write $z=|z|e^{i\theta}$ and $t=re^{i\theta}$ with $r\in [0,|z|]$, for some
$\theta \in [0,2\pi]$. Now,
\begin{align*}
\left|\int_{0}^{z}g'(z-t)f'(t)dt\right| & = \left|\int_{0}^{|z|}g'((|z|-r)e^{i\theta})
f'(re^{i\theta})e^{i\theta}dr\right| \\
& \leq \int_{0}^{|z|}|g'((|z|-r)e^{i\theta})|\,|f'(re^{i\theta})|dr \\
& \leq \int_{0}^{|z|}\frac{\|f\|_{\mathcal{B}}}{1-(|z|-r)^2} \frac{\|g\|_{\mathcal{B}}}{1-r^2}dr \\
& \leq \int_{0}^{|z|}\frac{\|f\|_{\mathcal{B}}}{(1-|z|+r)} \frac{\|g\|_{\mathcal{B}}}{(1-r)}dr
\\ & = \|f\|_{\mathcal{B}}\|g\|_{\mathcal{B}} \frac{2}{2-|z|}\log \frac{1}{1-|z|}.
\end{align*}
Since $(1-|z|)\log \frac{1}{1-|z|} \to 0$ as $|z|\to1$, 
 $(1-|z|)\log \frac{1}{1-|z|}$ is
 bounded on $\D$. Further, $\frac{2(1+|z|)}{2-|z|}\leq 4$ for all $z\in\D$. Hence,
\begin{align}
\label{fgderivduhaeq}(1-|z|^2)\left|\int_{0}^{z}g'(z-t)f'(t)dt\right| &
\leq \|f\|_{\mathcal{B}}\|g\|_{\mathcal{B}} \frac{2(1+|z|)}{2-|z|}(1-|z|)\log \frac{1}{1-|z|} \\
& \leq C\|f\|_{\mathcal{B}}\|g\|_{\mathcal{B}} \ \ \ \text{for some}~ C >0.\nonumber
\end{align}
Now \eqref{eqduhderiv} implies,
\begin{align}
(1-|z|^2)|(f\circledast g)'(z)|  & \leq (1-|z|^2)\left|\int_{0}^{z}g'(z-t)f'(t)dt\right| \nonumber \\
\label{fduhgderiv}& \ \ +|g(0)|(1-|z|^2)|f'(z)|+|f(0)|(1-|z|^2)|g'(z)| \\
& \leq C\|f\|_{\mathcal{B}}\|g\|_{\mathcal{B}}+ \|g\|_{\mathcal{B}}\|f\|_{\mathcal{B}} +
 \|f\|_{\mathcal{B}}\|g\|_{\mathcal{B}} < \infty. \nonumber 
\end{align}
It shows, $f\circledast g\in \mathcal{B}$.
Other properties of $\mathcal{B}$ as a Banach algebra under the Duhamel product are obvious
to verify and we are done.
 \epf

 If $f$ and $g$ are in $\mathcal{B}_0$, then equations \eqref{fgderivduhaeq} and \eqref{fduhgderiv}
 implies that $f\circledast g$ is also in $\mathcal{B}_0$ and hence we have the following corollary:
 \bcor\label{Boduhalg}
The little Bloch space, $\mathcal{B}_0$ is a Banach subalgebra of Duhamel algebra $\mathcal{B}$.
 \ecor

We end with the following concluding remark. We have characterized when 
a composition operator is Duhamel multiplicative. It is interesting to know
what are all operators which are Duhamel multiplicative. More specifically,
Characterize all bounded operators $T: H^p\to H^p$ which are  multiplicative
under Duhamel product.

\end{document}